\documentclass[12pt]{amsart} 
\usepackage{amsmath, amsfonts, amssymb, latexsym, amscd, enumerate}
\usepackage{pictexwd,dcpic}
\usepackage{graphicx}

\newtheorem{theorem}{Theorem}[section]
\newtheorem{lemma}[theorem]{Lemma}
\newtheorem{proposition}[theorem]{Proposition}
\newtheorem{corollary}[theorem]{Corollary}
\newtheorem{remark}[theorem]{Remark}
\newtheorem{example}[theorem]{Example}
\newtheorem{definition}[theorem]{Definition}
\newtheorem{notation}[theorem]{Notation}

\numberwithin{equation}{section}

\def\zN{\mathbb N}
\def\zZ{\mathbb Z}
\def\zQ{\mathbb Q}
\def\zR{\mathbb R}
\def\zC{\mathbb C}
\def\zT{\mathbb T}
\def\cA{\mathcal A}
\def\cE{\mathcal E}
\def\cF{\mathcal F}
\def\cL{\mathcal L}
\def\cK{\mathcal K}
\def\cI{\mathcal I}
\def\cO{\mathcal O}
\def\cT{\mathcal T}

\def\fG{\mathfrak G}
\def\fI{\mathfrak I}
\def\sq{$\text{}\hfill\square$}

\def\pf{\it Proof. \rm}
\def\st{\,\bigl\vert\,}

\newcommand{\abs}[1]{\vert#1\vert}
\newcommand{\norm}[1]{\vert\vert#1\vert\vert}
\newcommand{\innprod}[2]{\langle#1, #2\rangle}
\newcommand{\subgp}[1]{\langle#1\rangle}

\newcommand{\diag}[1]{\mbox{Diag}(#1)}

\title[Simplicity Criterion]{Simplicity criterion for $C^*$-algebras associated with topological group quivers}
\author{Shawn J. $\rm M^\MakeLowercase{c}$Cann}
\address{Department of Mathematics and Statistics, University of Regina, Regina, Sk, S4S 0A2}
\email{mccann1s@uregina.ca}
\begin{document}
\maketitle

\begin{abstract} Topological quivers generalize the notion of directed graphs in which
the sets of vertices and edges are locally compact (second countable) Hausdorff spaces. 
Associated to a topological quiver $Q$ is a $C^*$-correspondence, and in turn, a Cuntz-Pimsner 
algebra $C^*(Q).$ Given $\Gamma$ a locally compact group and $\alpha$ and $\beta$ 
endomorphisms on $\Gamma,$ one may construct a topological quiver $Q_{\alpha,\beta}(\Gamma)$
with vertex set $\Gamma,$ and edge set $\Omega_{\alpha,\beta}(\Gamma)=
\{(x,y)\in\Gamma\times\Gamma\st \alpha(y)=\beta(x)\}.$ In \cite{Mc1}, the author
examined the Cuntz-Pimsner algebra $\cO_{\alpha,\beta}(\Gamma):=C^*(Q_{\alpha,\beta}(\Gamma))$
and found  generators (and their relations) of $\cO_{\alpha,\beta}(\Gamma).$ In this paper, the author translates a known criterion for simplicity of topological quivers into a precise criterion for the simplicity of topological group relations.
\end{abstract}

\section{Introduction and Notation}
	\subsection{Introduction}In 2005, Muhly and Tomforde \cite{MT} defined a generalization of a directed graph called a \it topological quiver\rm. This is a
5-tuple $Q=(X,E,r,s,\lambda)$ where $X$ and $E$ are locally compact (second countable) Hausdorff
spaces, $r$ and $s$ are continuous maps from $X$ to $E$ with $r$ open, and $\lambda = \{\lambda_x\}_{x\in E}$
is a system of Radon measures. One can then create a corresponding Cuntz-Pimsner $C^*$-algebra $C^*(Q).$
One rather important property of a $C^*$-algebra is that of simplicity and Muhly and Tomforde calculate that 
$C^*(Q)$ is simple if and only if all open subsets $U\subseteq X$ with the property $s(e)\in E$ implies $r(e)\in U$
must be either $X$ or the empty set and if the set of all base points of loops in $Q$ with no exit has empty interior.
In general this can be a tedious and often difficult procedure to check.

In \cite{McThesis,Mc1,Mc2}, the author defines a particularly interesting topological quiver
$$Q=(\Gamma,\Omega_{\alpha,\beta}(\Gamma),r,s,\lambda)$$ 
where $\Gamma$ is a locally compact group, $\alpha$ and $\beta$ are endormorphism of $\Gamma,$
$$\Omega_{\alpha,\beta}(\Gamma)=\{(x,y)\in\Gamma\times\Gamma\st \alpha(y)=\beta(x)\}$$
and $\lambda$ is an appropriate family of Radon measures. Then forms the relative Cuntz-Pimsner Algebra denoted
$$\cO_{\alpha,\beta}(\Gamma).$$
In \cite{Mc1}, the spatial structure including generators, relations, and co-limit structure are examined. In \cite{Mc2},
a six-term exact sequence is prepared for the $K$-groups of $\cO_{\alpha,\beta}(\Gamma)$ and is used to calculate
the $K$-groups for $\cO_{F,G}(\zT^d),$ where $\zT^d$ is the $d$-torus and $F$ and $G$ are integral matrices with non-zero
determinant. With all this structure in place, the $C^*$-algebra is a classifiable (by its $K$-groups) Kirchberg algebra is 
it is also simple.

The aim of this paper is to recast the Muhly and Tomforde conditions into an easier condition for the $C^*$-algebras 
described above and furthermore, to determine the simplicity of $\cO_{F,G}(\zT^d)$ solely by the properties of $F$ and $G.$
That is, it will be shown that one merely needs to check the density of the connected component of the identity.
	\subsection{Notation}The sets of natural numbers, integers, rationals numbers, real numbers and complex numbers will be denoted by
$\zN$, $\zZ$, $\zQ$, $\zR$, and $\zC,$ respectively. Convention: $\zN$ does not contain zero.
$\zZ_0^+$ will denote the set $\zN\cup\{0\},$ $\zR^+$ denotes the set $\{r\in\zR\st r>0\}$ and $\zR_0^+=\zR^+\cup\{0\}.$
Finally, $\zZ_p$ denotes the abelian group $\zZ/p\zZ=\{0,1,...,p-1\mod p\}$ and 
$\zT$ denotes the torus $\{z\in\zC\st \abs{z}=1\}.$ Whenever convenient, view $\zZ_p\subset\zT$ by 
$\zZ_p\cong\{z\in\zT\st z^p=1\}.$

For a topological space $Y$, the closure of $Y$ is denoted $\overline{Y}.$ Given a locally compact Hausdorff space $X$, let
\begin{enumerate} 
\item $C(X)$ be the continuous complex functions on $X$;\index{C$(X)$}
\item $C_b(X)$ be the continuous and bounded complex functions on $X$;\index{C$\mbox{}_b(X)$}
\item $C_0(X)$ be the continuous complex functions on $X$ vanishing at infinity;\index{C$\mbox{}_0(X)$}
\item $C_c(X)$ be the continuous complex functions on $X$ with compact support.\index{C$\mbox{}_c(X)$}
\end{enumerate}
The supremum norm is denoted $\norm{\cdot}_\infty$ and defined by
$$\norm{f}_\infty=\sup_{x\in X}\{\abs{f(x)}\}$$
for each continuous map $f:X\to\zC.$ For a continuous function $f\in C_c(X),$ denote the open support of $f$ by 
$\mbox{osupp }f=\{x\in X\st f(x)\ne 0\}$ and the support of $f$ by $\mbox{supp }f=\overline{\mbox{osupp} f}.$  

For $C^*$-algebras $A$ and $B$, $A$ is isomorphic to $B$ will be written $A\cong B;$ for example, we use
$C(\zT^d)\otimes M_{N}(\zC)\cong M_{N}(C(\zT^d)).$ 
Moreover, $A^{\oplus n}$ denotes the $n$-fold direct sum $A\oplus\cdots\oplus A.$ 
Given a group $\Gamma$ and a ring $R$, a normal subgroup, $N$, of $\Gamma$ is denoted $N\lhd\Gamma$ and
an ideal, $I$, of $R$ is denoted $I\lhd R.$ Note if $R$ is a $C^*$-algebra then the term ideal denotes a closed two-sided ideal. 
Furthermore, $\mbox{End}(\Gamma)$ ($\mbox{End}(R)$) and $\mbox{Aut}(\Gamma)$
($\mbox{Aut}(R)$) denotes the set of endomorphisms of $\Gamma$ ($R$) and automorphisms of $\Gamma$ $(R$), respectively. 

Let $\alpha\in C(X)$ then $\alpha^\#\in\mbox{End(C(X))}$ denotes the endomorphism of $C(X)$ defined by
$$\alpha^\#(f)=f\circ\alpha\qquad\mbox{for each $f\in C(X)$}.$$ The set of $n$ by $n$ matrices with coefficients in a set $R$ will be denoted $M_n(R)$ and for any $F\in M_n(R),$
the transpose of $F$ is denoted $F^T$.
\section{Preliminairies}\subsection{Hilbert $C^*$-modules}Further details can be found in \cite{lan, RW}.

\begin{definition}\label{Hbmod}\cite{lan} \rm  If $A$ is a $C^*$-algebra, then a \emph{(right) Hilbert $A$-module}\index{Hilbert $C^*$-module} is a Banach space $\cE_A$
together with a right action of $A$ on $\cE_A$ and an $A$-valued inner product $\innprod{\cdot}{\cdot}_A$
satisfying
\begin{enumerate}
\item $\innprod{\xi}{\eta a}_A=\innprod{\xi}{\eta}_A a$
\item $\innprod{\xi}{\eta}_A =\innprod{\eta}{\xi}_A^*$
\item $\innprod{\xi}{\xi}\ge 0$ and $\norm{\xi}=\norm{\innprod{\xi}{\xi}_A^{1/2}}_A$
\end{enumerate}
for all $\xi$, $\eta\in\cE_A$ and $a\in A$ (if the context is clear, we denote $\cE_A$ simply by $\cE$). 
For Hilbert $A$-modules $\cE$ and $\cF$, call a function $T:\cE\to\cF$ \emph{adjointable}
\index{Hilbert $C^*$-module! Adjointable Operator, $\cL(\cE,\cF)$}
\index{Adjointable Operator, $\cL(\cE,\cF)$}
if there is a function $T^*:\cF\to\cE$ such that
$\innprod{T(\xi)}{\eta}_A=\innprod{\xi}{T^*(\eta)}_A$ for all $\xi\in\cE$ and $\eta\in\cF$.
Let $\cL(\cE,\cF)$ denote the set of adjointable ($A$-linear) operators from $\cE$ to $\cF$. 
If $\cE=\cF$, then $\cL(\cE):=\cL(\cE,\cE)$ is a $C^*$-algebra (see \cite{lan}.)
Let $\cK(\cE, \cF)$ denote the closed two-sided ideal of \emph{compact operators}
\index{Hilbert $C^*$-module! Compact Operators, $\cK(\cE,\cF)$} 
\index{Compact Operators, $\cK(\cE,\cF)$}
given by
$$\cK(\cE,\cF):=\overline{\mbox{span}}\{\theta_{\xi,\eta}^{\cE,\cF}\st\xi\in\cE,\,\eta\in\cF\}$$
where $$\theta_{\xi,\eta}^{\cE,\cF}(\zeta)=\xi\innprod{\eta}{\zeta}_A\qquad\mbox{for each $\zeta\in\cE$}.$$ 
Similarly, $\cK(\cE):=\cK(\cE,\cE)$ and $\theta_{\xi,\eta}^\cE$ (or $\theta_{\xi,\eta}$ if understood) denotes 
$\theta_{\xi,\eta}^{\cE,\cE}$.
For Hilbert $A$-module $\cE$, the linear span of $\{\innprod{\xi}{\eta}\st\xi,\eta\in\cE\}$, denoted $\innprod{\cE}{\cE}$,
once closed is a two-sided ideal of $A$. Note that $\cE\innprod{\cE}{\cE}$ is dense in $\cE$ (\cite{lan}).  
The Hilbert module $A_A$ refers to the Hilbert module $A$ over itself, where $\innprod{a}{b}=a^*b$ for all $a,b\in A$.
\end{definition}

\begin{definition}\label{ONB}\cite{KW} \rm A subset $\{u_i\}_{i\in\cI}\subset \cE$ is called a \emph{basis}\index{Basis}
provided the following reconstruction formula holds for all $\xi\in\cE:$
$$ \xi=\sum_{i\in\cI} u_i\cdot\innprod{u_i}{\xi}\qquad(\mbox{in }\cE,\norm{\cdot}.)$$
If $\innprod{u_i}{u_j}=\delta_i^j$ as well, call $\{u_i\}_{i\in\cI}$ an \emph{orthonormal basis}
\index{Basis! Orthonormal} of $\cE$.
\end{definition}  

\begin{definition}\label{Hbmorph}\cite{BB4, BB5} \rm If $A$ and $B$ are $C^*$-algebras, then an \emph{$A-B$ $C^*$-correspondence}\index{C$\mbox{}^*$-correspondence} $\cE$  is a right Hilbert 
$B$-module $\cE_B$ together with a left action of $A$ on $\cE$ given by a $*$-homomorphism $\phi_A:A\to\cL(\cE)$, 
$a\cdot\xi=\phi_A(a)\xi$ for $a\in A$ and $\xi\in\cE$. We may occasionally write, $_A\cE_B$ to denote an $A-B$ 
$C^*$-correspondence and $\phi$ instead of $\phi_A$.
Furthermore, if $_{A_1}\cE_{B_1}$ and $_{A_2}\cF_{B_2}$ are $C^*$-correspondences, then
a \emph{morphism}\index{C$\mbox{}^*$-correspondence! Morphism} 
$(\pi_1, T, \pi_2):\cE\to\cF$ consists of $*$-homomorphisms $\pi_i:A_i\to B_i$ and a linear map
$T:\cE\to\cF$ satisfying
\begin{enumerate}
\item[(i)] $\pi_2(\innprod{\xi}{\eta}_{A_2})=\innprod{T(\xi)}{T(\eta)}_{B_2}$
\item[(ii)] $T(\phi_{A_1}(a_1)\xi)=\phi_{B_1}(\pi_1(a_1))T(\xi)$
\item[(iii)] $T(\xi)\pi_2(a_2)=T(\xi a_2)$
\end{enumerate} for all $\xi,\eta\in\cE$ and $a_i\in A_i$.
\end{definition}

\begin{notation}\rm When $A=B$, we refer to $_A\cE_A$ as a $C^*$-correspondence over $A$. For $\cE$ a 
$C^*$-correspondence over $A$ and $\cF$ a $C^*$-correspondence over $B$, a morphism 
 $(\pi,T,\pi):\cE\to\cF$ will be denoted by $(T,\pi)$. 
\end{notation}

\begin{definition}\cite{MT} \rm If $\cF$ is the Hilbert module $ _CC_C$ where $C$ is a $C^*$-algebra with the inner product 
$\innprod{x}{y}_B=x^*y$ then call a morphism $(T,\pi):$ $_A\cE_B\to C$ of Hilbert
modules a \emph{representation}\index{C$\mbox{}^*$-correspondence! Representation} of $_A\cE_B$ into $C.$ 
\end{definition}

\begin{remark}\rm Note that a representation of $_A\cE_B$ need only satisfying $(i)$ and $(ii)$ 
of definition \ref{Hbmorph} as it was unnecessary to require (iii) (see \cite[Remark 2.7]{Mc1}).
\end{remark}

A morphism of Hilbert modules $(T,\pi):\cE\to\cF$ yields a $*$-homomorphism $\Psi_T:\cK(\cE)\to\cK(\cF)$ by
$$\Psi_T(\theta_{\xi,\eta}^\cE)=\theta_{T(\xi),T(\eta)}^\cF$$
for $\xi,\eta\in\cE$ and if $(S,\sigma):\mathcal D\to\cE$, and $(T,\pi):\cE\to\cF$ are morphisms of Hilbert modules then
$\Psi_T\circ\Psi_S=\Psi_{T\circ S}$. In the case where $\cF=B$ a $C^*$-algebra, we may first identify $\cK(B)$ as $B$,
and a representation $(T,\pi)$ of $\cE$ in a $C^*$-algebra $B$ yields a $*$-homomorphism $\Psi_T:\cK(\cE)\to B$ given
by $$\Psi_T(\theta_{\xi,\eta})=T(\xi)T(\eta)^*.$$

\begin{definition}\cite{MT} \rm For a $C^*$-correspondence $\cE$ over $A$, denote the ideal $\phi^{-1}(\cK(\cE))$ of $A$ by $J(\cE),$\index{J$(\cE)$} and let $J_\cE=J(\cE)\cap(\ker \phi)^\perp$\index{J$\mbox{}_\cE$} where $(\ker\phi)^\perp $ is the 
ideal $\{a\in A\st ab=0 \mbox{ for all }b\in\ker\phi\}$ .
If $_A\cE_A$ and $_B\cF_B$ are $C^*$-correspondences over $A$ and $B$ respectively and $K\lhd J(\cE)$, a morphism
$(T,\pi):\cE\to\cF$ is called \emph{coisometric on $K$}\index{C$\mbox{}^*$-correspondence! Representation! Coisometric on $K$} if $$\Psi_T(\phi_A(a))=\phi_B(\pi(a))$$
for all $a\in K$, or just \emph{coisometric},\index{C$\mbox{}^*$-correspondence! Representation! Coisometric} if $K=J(\cE)$.
\end{definition}

\begin{notation}\rm We denote $C^*(T,\pi)$ to be the $C^*$-algebra generated by $T(\cE)$ and $\pi(A)$ where
$(T,\pi):\cE\to B$ is a representation of $_A\cE_A$ in a $C^*$-algebra $B$. Furthermore, if $\rho:B\to C$ is a 
$*$-homomorphism of $C^*$-algebras, then $\rho\circ (T,\pi)$ denotes the representation $(\rho\circ T,\rho\circ\pi)$ of $\cE$.
\end{notation}

\begin{definition}\label{CPDefs}\cite{MT} \rm A morphism $(T_\cE,\pi_\cE)$ coisometric on an ideal $K$ is said to be \emph{universal}\index{C$\mbox{}^*$-correspondence! Representation! Universal} 
if whenever $(T,\pi):\cE\to B$ is a representation coisometric on $K$, there exists a $*$-homomorphism 
$\rho:C^*(T_\cE,\pi_\cE)\to B$ with $(T,\pi)=\rho\circ(T_\cE,\pi_\cE)$.  The universal $C^*$-algebra 
$C^*(T_\cE,\pi_\cE)$ is called the \emph{relative Cuntz-Pimsner algebra}\index{Cuntz-Pimsner Algebra}\index{$\cO(K,\cE)$}
\index{$\cO_\cE$} 
of $\cE$ determined by the ideal $K$ and 
denoted by $\cO(K,\cE)$. If $K=0$, then $\cO(K,\cE)$ is denoted
by $\cT(\cE)$ and called the \emph{universal Toeplitz $C^*$-algebra}\index{Toeplitz-Pimsner Algebra}\index{$\cT(\cE)$} for $\cE$. We denote $\cO(J_\cE,\cE)$ by $\cO_\cE$. 
\end{definition}
	\subsection{Topological Quivers}\begin{definition}\label{TQ}\cite{MT} \rm A \emph{topological quiver}\index{Topological Quiver} 
(or \emph{topological directed graph}\index{Topological Directed Graph}) $Q=(X,E,Y,r,s,\lambda)$ is a diagram
$$\begindc{\commdiag}[5]%50
\obj(10,0){$E$}
\obj(0,0){$X$}
\obj(20,0){$Y$}
\mor{$E$}{$X$}{$s$}[-1,0]
\mor{$E$}{$Y$}{$r$}
\enddc$$
where $X,E,$ and $Y$ are second countable locally compact Hausdorff spaces, $r$ and $s$ are continuous maps with $r$ open, along with a family $\lambda=\{\lambda_y\vert y\in Y\}$ of Radon measures on $E$ satisfying
\begin{enumerate}
\item $\mbox{supp }\lambda_y=r^{-1}(y)$ for all $y\in Y$, and
\item $y\mapsto\lambda_y(f)=\int_Ef(\alpha)d\lambda_y(\alpha)\in C_c(Y)$ for $f\in C_c(E).$
\end{enumerate}
\end{definition}

\begin{remark}\rm If $X=Y$ then write $Q=(X,E,r,s,\lambda)$ in lieu of $(X,E,X,r,s,\lambda).$
\end{remark}

\begin{remark}\rm The author provides a broad history and a series of examples of topological quivers in \cite{McThesis, Mc1}.
\end{remark}

Given a topological quiver $Q=(X,E,Y,r,s,\lambda)$, one may associate a correspondence  
$\cE_Q$ of the $C^*$-algebra $C_0(X)$ to the $C^*$-algebra $C_0(Y)$. Define left and right actions 
$$(a\cdot\xi\cdot b)(e)=a(s(e))\xi(e)b(r(e))$$
by $C_0(X)$ and $C_0(Y)$ respectively on $C_c(E)$. Furthermore, define the $C_c(Y)$-valued inner product
$$\innprod{\xi}{\eta}(y)=\int_{r^{-1}(y)}\overline{\xi(\alpha)}\eta(\alpha)d\lambda_y(\alpha)$$ 
for $\xi,\eta\in C_c(E)$, $y\in Y,$ and let $\cE_Q$\index{C$\mbox{}^*$-correspondence! Associated with a Topological Quiver} be the completion of $C_c(E)$ with respect to the norm 
$$\norm{\xi}=\norm{\innprod{\xi}{\xi}^{1/2}}_\infty=\norm{\lambda_y(\abs{\xi}^2)}_\infty^{1/2}.$$ 

\begin{definition}\rm Given topological quiver $Q$ over a space $X$, define the $C^*$-algebra, $C^*(Q)$
\index{C$\mbox{}^*$-algebra Associated with a Topological Quiver}\index{C$\mbox{}^*(Q)$}
\index{Topological Quiver! $C^*$-algebra Associated with} 
associated with $Q$ to be the Cuntz-Pimnser $C^*$-algebra $\cO_{\cE_Q}$\index{Cuntz-Pimsner Algebra} 
of the correspondence $\cE_Q$ over $A=C_0(X)$.
\end{definition}

	\subsection{Topological Group Quivers}\begin{definition}\cite{Mc1, McThesis}\label{TopGrpQuiver} \rm Let $\Gamma$ be a (second countable) locally compact group and let $\alpha,\beta\in\mbox{End}(\Gamma)$ be continuous. Define the closed subgroup, $\Omega_{\alpha,\beta}(\Gamma),$ 
of $\Gamma\times\Gamma,$\index{$\Omega_{\alpha,\beta}(\Gamma)$}
$$\Omega_{\alpha,\beta}(\Gamma)=\{(x,y)\in \Gamma\times\Gamma\st \alpha(y)=\beta(x)\}$$
and let $Q_{\alpha,\beta}(\Gamma)=(\Gamma,\Omega_{\alpha,\beta}(\Gamma), r,s,\lambda)$
\index{Q$\mbox{}_{\alpha,\beta}(\Gamma)$}\index{Topological Quiver}\index{Topological Quiver! $C^*$-algebra Associated with}  where $r$ and $s$ are the group homomorphisms defined by 
$$r(x,y)=x \qquad\mbox{and}\qquad s(x,y)=y$$ 
for each $(x,y)\in\Omega_{\alpha,\beta}(\Gamma)$ and $\lambda_x$ for $x\in\Gamma$ is the measure on 
$$r^{-1}(x)=\{x\}\times\alpha^{-1}(\beta(x))$$ 
defined by 
$$\lambda_x(B)=\mu(y^{-1}s(B\cap r^{-1}(x))\cap\ker\alpha)\qquad\mbox{(for any}\, y\in\alpha^{-1}(\beta(x)))$$
for each measurable $B\subseteq\Omega_{\alpha,\beta}(\Gamma)$ 
where $\mu$ is a left Haar measure (normalized if possible) on $r^{-1}(1_\Gamma)=\{1\}\times\ker\alpha$ 
(a closed normal subgroup of $\Gamma\times\Gamma;$ hence, a locally compact group). Note if $r^{-1}(x)=\emptyset$
then $\alpha^{-1}(\beta(x))=\emptyset$ and so $\lambda_x=0.$ 
This measure is well-defined, 
$$\mbox{supp }\lambda_x=\{x\}\times y\ker\alpha=\{x\}\times \alpha^{-1}(\beta(x))=r^{-1}(x)$$
and $y\mapsto \lambda_y(f)$ is a continuous compactly supported function (cf. \cite[Definition 3.1]{Mc1}.

Call $Q_{\alpha,\beta}(\Gamma)$ a \emph{topological group relation.}\index{Topological Group Relation}\index{Topological Relations}
Define $\cE_{\alpha,\beta}(\Gamma)$\index{$\cE_{\alpha,\beta}(\Gamma)$}\index{C$\mbox{}^*$-correspondence} 
to be the $C_0(\Gamma)$-correspondence $\cE_{Q_{\alpha,\beta}(\Gamma)}$ and
form the Cuntz-Pimsner algebra\index{$\cO_{\alpha,\beta}(\Gamma)$! Definition}\index{Cuntz-Pimsner Algebra}
$$\cO_{\alpha,\beta}(\Gamma):=C^*(Q_{\alpha,\beta}(\Gamma))=\cO(J_{\cE_{\alpha,\beta}(\Gamma)},\cE_{\alpha,\beta}(\Gamma))$$
and the Toeplitz-Pimsner algebra\index{$\cT_{\alpha,\beta}(\Gamma)$! Definition}\index{Toeplitz-Pimsner Algebra}
$$\cT_{\alpha,\beta}(\Gamma):=\cT(Q_{\alpha,\beta}(\Gamma)).$$
\end{definition}

\begin{remark}\rm It will be implicitly assumed that $\Gamma$ is second countable. Furthermore, since $\Gamma$ is locally compact Hausdorff, $r^{-1}(x)$ is closed and locally compact. Moreover, whenever $r$ is a local homeomorphism, $r^{-1}(x)$ is discrete and hence, $\lambda_x$ is counting measure (normalized when $\abs{\ker\alpha}<\infty$.)
\end{remark}\begin{example}[\cite{Mc1}]\label{ToriQuiver}\rm\index{Topological Group Relation! $\zT^d$-quivers}\index{Topological Group Relation}
\index{Topological Relations}\index{Topological Quiver}\index{$\cO_{F,G}(\zT^d)$! Definition}
For the compact abelian group $\zT^d,$ note $\mbox{End}(\zT^d)\cong M_d(\zZ)$ (\cite{W}); that is, 
an element $\sigma\in\mbox{End}(\zT^d)$ is of the form $\sigma_F$ for some $F\in M_d(\zZ)$ where 
$$\sigma_F(e^{2\pi it})=e^{2\pi i Ft}$$\index{$\sigma_F$}
for each $t\in\zZ^d.$ To simplify notation, use $F$ and $G$ in place of $\sigma_F$ and $\sigma_G$ whenever convenient.
For instance,
$$Q_{F,G}(\zT^d):=Q_{\sigma_{F},\sigma_{G}}(\zT^d)$$
and the $C^*$-correspondence
$$\cE_{F,G}(\zT^d):=\cE_{\sigma_F,\sigma_G}(\zT^d)$$
where $F,G\in M_d(\zZ)$. We will consider the cases when these maps are surjective; that is,
$\det F$ and $\det G$ are non-zero.

Let $F, G\in M_d(\zZ)$ where $\det F,\det G\ne 0$. Then $\abs{\ker\sigma_F}=\abs{\det F}$
and so, the $C(\zT^d)$-valued inner product becomes
$$\innprod{\xi}{\eta}(x)=\frac{1}{\abs{\det F}}\sum_{\sigma_{F}(y)=\sigma_{G}(x)}\overline{\xi(x,y)}\eta(x,y)$$
for $\xi,\eta\in\cE_{F,G}(\zT^d)$ and $x\in\zT^d.$
This is a finite sum since the number of solutions, $y,$ 
to $\sigma_F(y)=\sigma_G(x)$ given any $x\in\zT^d$ is $\abs{\det F}<\infty.$
The left action $\phi,$ defined by
$$\phi(a)\xi(x,y)=a(y)\xi(x,y)$$
for $a\in C(\zT^d),$ $\xi\in C(\Omega_{F,G}(\zT^d))$ and $(x,y)\in\Omega_{F,G}(\zT^d),$ is injective (cf. \cite[Remark 3.18]{Mc1}). 
\end{example}

\begin{remark}\rm It was shown in \cite[Corollary 3.20]{Mc1} that one may assume the matrix $F$ is positive diagonal. 
\end{remark}

Let $F=\diag{a_1,...,a_d}\in M_d(\zZ), G=(b_{jk})_{j,k=1}^d\in M_d(\zZ)$ where $a_j>0$ for each $j=1,...,d,$ $\det G\ne0$ and
let $G_j$ denote the $j$-th row of $G$, $(b_{jk})_{k=1}^d.$ Further, let $N=\det F=\prod_{j=1}^d a_j>0$ and let $$\fI(F)=\{\nu=(\nu_j)_{j=1}^d\in\zZ^d\st 0\le \nu_j\le a_j-1\}.\index{$\fI(F)$}$$ 
The $C(\zT^d)$-valued  inner product becomes
$$\innprod{\xi}{\eta}(x)=\frac{1}{N}\sum_{\sigma_F(y)=\sigma_G(x)} \overline{\xi(x,y)}\eta(x,y)$$
for all $\xi,\eta\in C(\Omega_{F,G}(\zT^d))$ and $x\in \zT^d.$

Given $\nu\in\fI(F)$, define $u_\nu\in C(\Omega_{F,G}(\zT^d))$ by 
$$u_\nu(x,y)=y^\nu=\prod_{j=1}^d y^{\nu_j}$$
for $(x,y)\in \Omega_{F,G}(\zT^d).$ It was shown in \cite{Mc1} that $\{u_\nu\}_{\nu\in\fI(F)}$ is a basis for $\cE_{F,G}(\zT^d)$
and also the following:

\begin{theorem}\cite[Theorem 3.23]{Mc1}\label{TdQuiverGen}\rm \index{$\cO_{F,G}(\zT^d)$}\index{$\cO_{\alpha,\beta}(\Gamma)$! Presentation}
\index{Cuntz-Pimsner Algebra}
Let $F=\diag{a_1,...,a_d}, G\in M_d(\zZ)$ where $\det F, \det G\ne0$ and let $G_j$ be the $j$-th row vector of $G$. Further, let $\fI(F)$ denote the set $\{\nu=(\nu_j)_{j=1}^d\in\zZ^d\st 0\le \nu_j\le a_j-1\}$. Then
$\cO_{F,G}(\zT^d)$ is the universal $C^*$-algebra generated by isometries $\{S_\nu\}_{\nu\in\fI(F)}$ and (full spectrum) commuting unitaries $\{U_j\}_{j=1}^d$ that satisfy the relations
\begin{enumerate}
\item $S_\nu^*S_{\nu^\prime}=\innprod{u_\nu}{u_{\nu^\prime}}=\delta_{\nu}^{\nu^\prime},$
\item $U^\nu S=S_\nu$ for all $\nu\in\fI(F),$
\item $U_j^{a_j}S=SU^{G_j},$ for all $j=1,...,d$ and
\item $1=\sum_{\nu\in\fI(F)} S_\nu S_\nu^*=\sum_{\nu\in\fI(F)} U^\nu SS^*U^{-\nu}$
\end{enumerate}
where $U^\nu$ denotes $\prod_{j=1}^dU_j^{\nu_j}.$ Furthermore, $\cT_{\alpha,\beta}(\Gamma)$ is the universal $C^*$-algebra generated by isometries $\{S_\nu\}_{\nu\in\fI(F)}$ and commuting unitaries $\{U_j\}_{j=1}^d$ that satisfy relations (1)-(3)
\end{theorem}
\section{Simplicity Criterion}We now recast known conditions for the simplicity of Cuntz-Pimsner algebras for the correspondences in our specific context.

\begin{definition}\cite{MT} \rm \index{Topological Quiver}
Let $Q=(X,E,r,s,\lambda)$ be a topological quiver. A \emph{path}\index{Path} in $Q$ is a finite sequence 
$\{e_1,...,e_n\}\subseteq E$ denoted
$e=e_1...e_n$ with $r(e_i)=s(e_{i+1})$ for all $i=1,...,n-1$ with \emph{length}\index{Length of Path} 
$\abs{e}=n$. Denote the set of paths of length $n$ by $E^n$.

If $e=e_1...e_n\in E^n$ and $r(e_n)=s(e_1)$ then call $e$ a \emph{loop}\index{Loop} of length $n$. The vertex $x=s(e_1)=r(e_n)$ is refered to as the \emph{base point}\index{Base Point} of the loop $e$. An \emph{exit}\index{Exit of a Loop}
of a loop $e=e_1...e_n$ is any $e_0\in E$ where for some $k\in\{1,...,n\},$ $s(e_0)=s(e_k)$ ($e_0\not= e_k.$) 
\end{definition}

\begin{definition}\cite{MT} \rm
Say $Q$ has \emph{condition ($L$)}\index{Condition $(L)$} if the set of base points of loops in $Q$ with no exits has empty interior.
\end{definition}

\begin{definition}\cite{MT} \rm  Let $Q=(X,E,r,s,\lambda)$ be a topological quiver. An open subset $U\subseteq X$ is 
\emph{hereditary}\index{Hereditary} if $s(e)\in U$ for some $e\in E$ guarantees that $r(e)\in U$. 
Say a hereditary subset $U\subseteq X$ is \emph{saturated}\index{Saturated} if 
$x\in X_{reg}$ and $r(s^{-1}(x))\subseteq U$ imply $x\in U$.
\end{definition} 

\begin{definition}\cite{MT} \rm A topological quiver, $Q=(X,E,r,s,\lambda)$, is \emph{minimal}\index{Minimal} 
if the only saturated hereditary open subsets of $Q$ are $X$ and $\emptyset.$
\end{definition}

\begin{theorem}\cite[Theorem 10.2]{MT} \rm \label{MuhlyThm} If $Q=(X,E,r,s,\lambda)$ is a topological quiver then $C^*(Q)$ is simple if and only if $Q$ is minimal and satisfies condition ($L$).
\end{theorem}

We now focus our attention to determining the necessary and sufficient conditions to ensure minimality and condition ($L$) for
the topological group quiver $Q_{\alpha,\beta}(\Gamma)$.
The following lemma may be found in \cite{W}.

\begin{lemma}\cite{W} \rm 
Every continuous surjective endomorphism of a compact group $\Gamma,$ $\alpha:\Gamma\to \Gamma,$ 
is measure-preserving.\\ \index{Measure Preserving}
\pf Let $m$ denote the normalized Haar measure of $\Gamma$ and define a probability measure $\mu(B)=m(\alpha^{-1}(B))$ for all Borel subsets $B$ of $\Gamma$. The measure $\mu$ is regular. 
We have
$$\mu(\alpha(x)\cdot B)=m(\alpha^{-1}(\alpha(x)\cdot B))=m(x\cdot\alpha^{-1}(B))=m(\alpha^{-1}(B))=\mu(B)$$
for all Borel subsets $B$ of $\Gamma$. Since $\alpha$ maps $\Gamma$ onto $\Gamma$, we see that $\mu$ is left translation invariant and thus, by uniqueness of the Haar measure, $\mu=m$; that is, $m(\alpha^{-1}(B))=\mu(B)=m(B)$ for all Borel subset $B$ of $\Gamma$.\\\sq
\end{lemma}   

\begin{theorem} \rm\label{L} Let $\alpha,\beta\in\mbox{End}(\Gamma)$. If either
\begin{enumerate}
\item $\beta$ is not an injection, or
\item $\alpha$ and $\beta$ are surjections, $\alpha$ is not an injection and $\abs{\ker\alpha}<\infty$
\end{enumerate} 
then $Q_{\alpha,\beta}(\Gamma)$ satisfies condition ($L$).\\ \index{Condition $(L)$}
\pf If $e$ is a loop of length $n$ in $Q_{\alpha,\beta}(\Gamma)$, then $e$ has the following form:
$$e=e_1...e_n$$ 
where 
$$e_k=(x_k,x_{k-1})\mbox{ for }k=1,...,n,$$ 
$$r(e_n)=x_n=x_0=s(e_1)$$ 
and $$x_k=r(e_k)=s(e_{k+1})\in\alpha^{-1}(\beta(x_{k+1}))\mbox{ for }k=1,...,n-1.$$
But if $\beta$ is not an injection then there is a non-trivial $z\in\ker\beta$. Furthermore, for any $k\in\{1,...,n\}$,
let $f_k=(zx_{k}, x_{k-1}).$ Then $\beta(zx_k)=\beta(x_k)=\alpha(x_{k-1})$ and so, $f_k\in\Omega_{\alpha,\beta}(\Gamma)$ with the property that $s(f_k)=s(e_k)$ and $f_k\not= e_k.$ Hence, $e$ has an exit and as $e$ was arbitrary, all loops have exits. Thus,
the set of base points of loops without exits is empty and hence, has empty interior.

If, alternatively, $\beta$ is an automorphism (satisfies (2) and not (1)), 
then there exists an inverse of $\beta$, call it $\beta^{-1}$ and so, define
$$B_n=\{x\in \Gamma\st (\beta^{-1}\alpha)^n(x)=x\}\qquad(\mbox{base points of loops of length $n$})$$
for each $n\in\zN.$ It is now evident that the set of base points of loops, $L,$ is 
$$L=\bigcup_{n\in\zN}B_n.$$
Furthermore, note $B_n$ is closed since $(\beta^{-1}\alpha)^n$ is continuous, and hence $B_n$ is measurable. Let 
$\alpha_0=\beta^{-1}\alpha\in\mbox{End}(G)$ for the ease of this argument. Note $\alpha_0$ is a continuous surjective endomorphism, hence by the above lemma, $\alpha_0$ is measure-preserving. Furthermore,
$$\alpha_0^{-1}(B_n)=\ker\alpha\cdot B_n.$$
Indeed, if $\alpha_0(z)\in B_n$ then $\alpha_0^{n+1}(z)=\alpha_0(z)$ and hence, $\alpha_0^n(z)\in\alpha_0^{-1}(\alpha_0(z))=(\ker\alpha_0)z=(\ker\alpha)z.$ So $\alpha_0(z)=z_0z$ for some $z_0\in\ker\alpha$ and thus,
$\alpha_0^n(z_0z)=z_0z$ implies $z_0z\in B_n$ and hence, $z\in(\ker\alpha)\cdot B_n.$ Conversely, if $z_0z\in B_n$ for some
$z_0\in\ker\alpha$ then $\alpha_0^{n+1}(z)=\alpha_0(\alpha_0^n(z_0z))=\alpha_0(z_0z)=\alpha_0(z)$ and thus, 
$\alpha_0(z)\in B_n.$

Also, for $z_1,z_2\in\ker\alpha$ with $z_1\not= z_2$, and $x,y\in B_n$, $z_1x=z_2y$ implies $$x=\alpha_0^n(x)=\alpha_0^n(z_1x)=\alpha_0^n(z_2y)=\alpha_0^n(y)=y$$
and thus, $z_1=z_2yx^{-1}=z_2;$ that is, it is shown that $z_1\cdot B_n$ is disjoint from $z_2\cdot B_n$ whenever $z_1\not= z_2,$ for $z_1,z_2\in\ker\alpha.$     

Hence, for Haar measure $m$,
$$m(B_n)=m(\alpha_0^{-1}(B_n))=m(\ker\alpha\cdot B_n)=m(\bigsqcup_{z\in\ker\alpha} z\cdot B_n)=\sum_{z\in\ker\alpha}m(B_n)
=\abs{\ker\alpha}\,m(B_n).$$
Since $1<\abs{\ker\alpha}<\infty$, it must be that $m(B_n)=0;$ that is, $B_n$ must have empty interior. Moreover, since $\Gamma$ is a Baire space, $L=\cup_n B_n$ has empty interior. Therefore, $Q_{\alpha,\beta}(\Gamma)$ satisfies condition $(L)$.\\\sq
\end{theorem}

\begin{remark}\rm (1) of Theorem \ref{L} does not require $\Gamma$ compact.
\end{remark}

%We shall now focus on determining when minimality is achieved, but to do so we need the following definition of an adjoint.

%\begin{definition}\rm Let $Q=(X,E,Y,r,s,\lambda)$ be a topological quiver such that $s$ is also an open map. Then an adjoint quiver %of $Q$, denote $Q^*$, is a topological quiver $Q^*=(Y,E,X,s,r,\lambda^\prime).$ Provided the existence of an adjoint quiver, we %call $Q$ adjointable.
%\end{definition}

%\begin{remark}\rm It should be stated that the adjoint of $Q$ is conditional on the existence of a family of Radon measures %$\lambda^\prime$ satisfying the topological quiver properties. Furthermore, the adjoint is none other than reversing the arrows.
%\end{remark}

%\begin{example}\rm For a topological group quiver, $Q_{\alpha,\beta}(G),$ $Q_{\beta,\alpha}(G)$ is an adjoint. Indeed,
%if $(x,y)\in\Omega_{\alpha,\beta}(G)$ then $(y,x)\in\Omega_{\beta,\alpha}(G)$ and so, we may choose $\lambda_x^\prime$
%to be (normalized) Haar measure on $s^{-1}(x)$ remarking that $s$ is an open map.
%\end{example}

Now, with the standing assumption that $\alpha$ and $\beta$ are surjective endomorphisms on a compact group $\Gamma$ with  $\abs{\ker\alpha},\abs{\ker\beta}<\infty$, we intend to characterize when $Q_{\alpha,\beta}(\Gamma)$ is minimal. But first, a lemma is needed.

\begin{lemma}\rm \label{Crit1} Let $\alpha$ and $\beta$ be surjective endomorphisms on a compact group $\Gamma$ with finite
kernels. An open set $U\subseteq \Gamma$ is a saturated hereditary subset of $Q_{\alpha,\beta}(\Gamma)$ if and only if $\alpha^{-1}(\beta(U))=\beta^{-1}(\alpha(U))=U.$\\\index{Saturated}\index{Hereditary}
\pf Let $U$ be a saturated hereditary open subset of $Q$, then certainly, 
$$r(s^{-1}(y))\subseteq U\mbox{ if and only if } y\in U.$$
Now realize $s^{-1}(y)=\{(x,y)\st \alpha(y)=\beta(x)\}$ and so, $r(s^{-1}(y))=\beta^{-1}(\alpha(y))$.
Hence, recharacterize $U$ by $\beta^{-1}(\alpha(y))\subseteq U$ if and only if $y\in U.$
Furthermore, if $u^\prime\notin U$ then, by the surjectivity of $\alpha$, there exists $y\in\Gamma$ such that $\alpha(y)=\beta(u^\prime).$ But note that $y\notin U$ since otherwise, if $y\in U$ then $\beta^{-1}(\alpha(y))\cap U^c\not=\emptyset.$ Hence, $$U^c\subseteq \beta^{-1}(\alpha(U^c)).$$
Next, if $z\notin\alpha(U)$ then  $\alpha^{-1}(z)\subseteq U^c$ and thus, $z\in\alpha(U^c)$. Conversely,
first note $\beta^{-1}(\alpha(\ker\alpha\cdot U))=\beta^{-1}(\alpha(U))\subseteq U,$ so $\ker\alpha\cdot U=U$ 
(remember $1_\Gamma\in\ker\alpha$.) Now if $\alpha(u)=\alpha(u^\prime)$ with $u\in U$, then 
$$u^\prime\in\alpha^{-1}(\alpha(U))=\ker\alpha\cdot U=U.$$
Hence, $\alpha(U)\cap\alpha(U^c)=\emptyset,$ but since $\alpha$ is surjective, $\alpha(U)\cup\alpha(U^c)=\Gamma$ and so,
$\alpha(U^c)=\alpha(U)^c$.

This now yields
$$U^c\subseteq\ker\alpha\cdot U^c=\alpha^{-1}(\alpha(U^c))=\alpha^{-1}((\alpha(U))^c)\subseteq U^c$$
and thus, $\alpha^{-1}(\alpha(U^c))=U^c$. Use the normalized Haar measure $\mu:$
$$\mu(U^c)=\mu(\alpha^{-1}(\alpha(U^c)))=\mu(\beta^{-1}(\alpha(U^c)))$$
and since both $U^c$ and $\beta^{-1}(\alpha(U^c))$ are closed with $U^c\subseteq \beta^{-1}(\alpha(U^c))$, it must be that
$\beta^{-1}(\alpha(U^c))=U^c$. Hence,
$$U=(\beta^{-1}(\alpha(U^c)))^c\subseteq \beta^{-1}(\alpha(U))\subseteq U;$$
that is, $\beta^{-1}(\alpha(U))=U.$ 

It should now be noted that it has been shown that $\beta(U)=\alpha(U)$, so $U\subseteq\alpha^{-1}(\beta(U))$.
If $\alpha(y)=\beta(u)$ for some $u\in U$, then there exists $u^\prime\in U$ such that $\beta(u)=\alpha(u^\prime)$.
hence, $y(u^\prime)^{-1}\in\ker\alpha$ and so $y=y(u^\prime)^{-1}u^\prime\in\ker\alpha\cdot U=U.$
Therefore, $\alpha^{-1}(\beta(U))=\beta^{-1}(\alpha(U))=U$.

As for the converse, $\beta^{-1}(\alpha(U))=U$, so it is enough to show that $\beta^{-1}(\alpha(y))\subseteq U$ implies $y\in U.$ To this end, assume $\beta^{-1}(\alpha(y))\subseteq U;$ that is, $\alpha(y)\in\beta(U)$ and so 
$$y\in\alpha^{-1}(\beta(U))=U.$$\sq
\end{lemma}

\begin{remark}\rm Note two properties for open saturated hereditary set $U$:
\begin{enumerate}
\item $\ker\alpha\cdot U=U$, and
\item $\ker\beta\cdot U=U$.
\end{enumerate}
(1) was proven in the previous lemma, but to see (2) note $\beta(\ker\beta\cdot U)=\beta(U)=\alpha(U)$, so by the previous lemma, $\ker\beta\cdot U\subseteq U$. Hence, $\ker\beta\cdot U=U.$ 
\end{remark}

\begin{definition}\label{GDef} \rm Let $\fG_0=\{1_{\Gamma}\},$ the trivial subgroup of $\Gamma,$ and let 
$$\fG_n=\begin{cases} \beta^{-1}(\alpha(\fG_{n-1}))\mbox{ if $n>0$}\\
\alpha^{-1}(\beta(\fG_{n+1}))\mbox{ if $n<0$}.
\end{cases}$$
Then let $\fG=\cup_{n\in\zZ}\fG_n$ and let $\Gamma_0$ be the subgroup (of $\Gamma$) generated by $\fG$.
\index{$\fG$}\index{$\Gamma_0$}
\end{definition}

We describe a few properties of $\fG$ and $\Gamma_0$.

\begin{lemma}\label{fGn}\rm The following are true for $\fG$ and $\Gamma_0$ defined above:
\begin{enumerate}
\item $\fG_n$ is a normal subgroup of $\Gamma$ for each $n\in\zZ$
\item $\fG_n\subseteq \fG_{n+1}$ for $n\ge 0$ and $\fG_n\subseteq \fG_{n-1}$ for $n\le0$.
\item $\Gamma_0$ is the set of finite products of elements in $\fG$
\item $\Gamma_0$ is a normal subgroup of $\Gamma$.
\end{enumerate}
\pf (1) is clear, (3) follows from (1) and (4) is evident by (1) and (3).

(2) Notice that $\{1_\Gamma\}=\fG_0\subseteq \fG_1=\ker\beta$. Proceed with induction. Assume that $\fG_{n-1}\subseteq\fG_n$ for $n>0,$ then given any $x\in \fG_n,$ there exists $y\in \fG_{n-1}\subseteq \fG_n$ such that $\beta(x)=\alpha(y)$; that is, $x\in\beta^{-1}(\alpha(\fG_n))=\fG_{n+1}.$ A similar argument for $n\le0$ shows that $\fG_n\subseteq \fG_{n-1}$.  \\\sq
\end{lemma}

Define the directed set $\zZ_0^+\times\zZ_0^+$ by 
$$(n_1,m_1)\le (n_2,m_2)\iff n_1\le n_2\mbox{ and } m_1\le m_2$$
for $(n_i,m_i)\in\zZ_0^+\times\zZ_0^+$ $(i=1,2.)$ Let 
$$\cF_{(n,m)}=\innprod{\fG_n}{\fG_{-m}},$$
the subgroup of $\Gamma$ generated by $\fG_n$ and $\fG_{-m}.$
For $d_i=(n_i,m_i)\in\zZ_0^+\times\zZ_0^+$ for $i=1,2$ where $(n_1,m_1)\le (n_2,m_2)$ in $\zZ_0^+\times\zZ_0^+,$ 
define the morphism
$$f_{d_1,d_2}:\cF_{d_1}\hookrightarrow\cF_{d_2},$$
the inclusion of $\cF_{d_1}$ into $\cF_{d_2}.$ One can verify that 
$$\mathcal C=(\{\cF_d\}_{d\in\zZ_0^+\times\zZ_0^+}, \{f_{d_1,d_2}\}_{d_1\le d_2})$$
defines a small category with objects $\mbox{Ob}(\mathcal C)=\{\cF_d\}_{d\in\zZ_0^+\times\zZ_0^+}$ and morphisms
$\mbox{Hom}(\mathcal C)
=\{f_{d_1,d_2}:\cF_{d_1}\hookrightarrow\cF_{d_2}\st d_1\le d_2, d_1,d_2\in\zZ_0^+\times\zZ_0^+\}.$

\begin{proposition}\label{SimpleColimit}\rm 
With the above notation, $\Gamma_0$ is the colimit (in the category of groups) of the diagram 
$\cF:\zZ_0^+\times\zZ_0^+\to\mathcal C$ defined by
$$\cF(n,m)=\cF_{(n,m)}=\innprod{\fG_n}{\fG_{-m}}$$
for $(n,m)\in\zZ_0^+\times\zZ_0^+$ and
$$\cF(d_1, d_2)=f_{d_1,d_2}:\cF_{d_1}\hookrightarrow\cF_{d_2}$$
for $d_1\le d_2$ in $\zZ_0^+\times\zZ_0^+.$\\
\pf For $d\in\zZ_0^+\times\zZ_0^+,$ let $\phi_d:\cF_d\hookrightarrow\Gamma_0$ be the inclusion of $\cF_d$ into $\Gamma_0.$
It is easy to see that
$$\begindc{0}[10]%100
\obj(4,5)[A]{$\cF_{d_1}$}
\obj(10,5)[B]{$\cF_{d_2}$}
\obj(7,0)[C]{$\Gamma_0$}
\mor{A}{B}{$f_{d_1,d_2}$}
\mor{A}{C}{$\phi_{d_1}$}[-1,0]
\mor{B}{C}{$\phi_{d_2}$}
\enddc$$
is a commutative diagram for any $d_1\le d_2$ in $\zZ_0^+\times\zZ_0^+.$
Suppose $(\cA, \psi)$ is co-cone of diagram $\cF$; that is,
$$\begindc{0}[10]%100
\obj(4,5)[A]{$\cF_{d_1}$}
\obj(10,5)[B]{$\cF_{d_2}$}
\obj(7,0)[C]{$\cA$}
\mor{A}{B}{$f_{d_1,d_2}$}
\mor{A}{C}{$\psi_{d_1}$}[-1,0]
\mor{B}{C}{$\psi_{d_2}$}
\enddc$$
is a commutative diagram for any $d_1\le d_2$ in $\zZ_0^+\times\zZ_0^+.$ Then for $g\in\Gamma_0,$ there exists 
$d=(n,m)\in\zZ_0^+\times\zZ_0^+$ such that $g\in\cF_d$ and so, let
$$\phi:\Gamma_0\to\cA$$
where
$$g\mapsto\psi_d(g).$$
This map is a well defined homomorphism. Let $d_i=(n_i,m_i)\in\zZ_0^+\times\zZ_0^+$ for $i=1,2$ be such that $g\in\cF_{d_1}\cap\cF_{d_2}.$ Then for $n=\max(n_1,n_2)$ and $m=\max(m_1,m_2),$ 
$$g\in\cF_{(n,m)}.$$
Note $d_i\le d=(n,m)$ for $i=1,2$ and so, 
$$\psi_{d_1}(g)=\psi_d(g)=\psi_{d_2}(g).$$ 
Furthermore, if $g,g^\prime\in\Gamma_0$ then there exists $d\in\zZ_0^+\times\zZ_0^+$ such that $g,g^\prime\in\cF_d$
and so, 
$$\phi(gg^\prime)=\psi_d(gg^\prime)=\psi_d(g)\psi_d(g^\prime)=\phi(g)\phi(g^\prime).$$
Note for any $d\in\zZ_0^+\times\zZ_0^+$ and $g\in\cF_d,$
$$(\phi\circ\phi_d)(g)=\phi(g)=\psi_d(g);$$
that is, 
$$\begindc{\commdiag}[5]%30
\obj(10,30)[A]{$\cF_{d_1}$}
\obj(40,30)[B]{$\cF_{d_2}$}
\obj(25,15)[C]{$\Gamma_0$}
\obj(25,0)[D]{$\cA$}
\mor{A}{B}{$f_{d_1,d_2}$}
\mor{A}{C}{$\phi_{d_1}$}[-1,0]
\mor{B}{C}{$\phi_{d_2}$}
\cmor((9,27)(11,8)(22,1))
\pright(5,15){$\psi_{d_1}$}
\cmor((41,27)(39,8)(28,1))
\pleft(45,15){$\psi_{d_2}$}
\mor{C}{D}{$\phi$}[-1,1]
\enddc$$
is a commutative diagram for any $d_1\le d_2\in\zZ_0^+\times\zZ_0^+.$
Finally note if $\phi^\prime:\Gamma_0\to\cA$ is a homomorphism making the above diagram commute, then given $g\in\Gamma_0$
and $d\in\zZ_0^+\times\zZ_0^+$ such that $g\in\cF_d,$
$$\phi^\prime(g)=(\phi^\prime\circ\phi_d)(g)=\psi_d(g)=\phi(g).$$
Hence, the homomorphism $\phi:\Gamma_0\to\cA$ is unique.\\\sq
\end{proposition}

\begin{proposition}\rm Let $U\in\Gamma$ be a saturated hereditary subset for the topological quiver $Q_{\alpha,\beta}(\Gamma)$. 
Then given $u\in U$ and $x\in\fG$,\index{Saturated}\index{Hereditary}
$ux,xu\in U.$\\
\pf Proceed by induction. Let $u\in U$ and $x\in \fG_n$ with $n\ge 0$. Then if $n=0$, clearly $xu=ux=u\in U$. Now assume that
$ux,xu\in U$ for every $x\in \fG_n$ and let $y\in \fG_{n+1}$. There exists $x\in\fG_n$ such that $\beta(y)=\alpha(x)$ and since
$\beta(U)=\alpha(U)$ (by Lemma \ref{Crit1}),  there exists $u^\prime\in U$ such that $\beta(u)=\alpha(u^\prime).$ Then
$$\beta(uy)=\alpha(u^\prime x)\qquad\mbox{and}\qquad \beta(yu)=\alpha(xu^\prime)$$
and by the inductive hypothesis, $u^\prime x,xu^\prime\in U$. So $uy,yu\in\beta^{-1}(\alpha(U))=U.$ A similar argument for $n\le0$ shows that $xu,ux\in U$ for $x\in \fG_n$.\\\sq
\end{proposition}

\begin{corollary}\rm For $U\subseteq\Gamma$ a saturated hereditary (open) subset for the topological quiver $Q_{\alpha,\beta}(\Gamma),$ 
$$U\cdot \Gamma_0=\Gamma_0\cdot U=U.$$
\pf Combine Lemma \ref{fGn} and the previous proposition.\\\sq
\end{corollary}

Recall that a second countable compact Hausdorff space is metrizable (see \cite{M}). We remark that if $D$ is a dense
subset of $\Gamma$ then so is any left (or right) translate of $D$ since left (or right) multiplication by $x\in\Gamma$ is a 
homeomorphism of $\Gamma.$

\begin{theorem}\rm \label{MinThm} $\Gamma_0$ is dense in $\Gamma$ if and only if $Q_{\alpha,\beta}(\Gamma)$ is minimal.\\
\index{Minimal}\index{Q$\mbox{}_{\alpha,\beta}(\Gamma)$! Minimal}
\pf First, assume $\Gamma_0$ is dense in $\Gamma.$ Let $U$ be a non-empty saturated hereditary open subset of $Q_{\alpha,\beta}(\Gamma).$ Since $\Gamma_0$ is dense, so is $x\cdot \Gamma_0$ for any $x\in \Gamma_0$. Hence, since $U$ is open, given $x\in \Gamma$ there exists $z\in\Gamma_0$ such that $x\cdot z\in U$. But $z^{-1}\in \Gamma_0$ implies $x=x\cdot z\cdot z^{-1}\in U$. Therefore, $U=\Gamma.$

Conversely, if $\Gamma_0$ is not dense in $\Gamma$, then let $U=\overline{\Gamma_0}^c$, an open non-empty subset of 
$\Gamma$ since $1_\Gamma\in\Gamma_0.$ We show that $U$ is a saturated hereditary subset of $\Gamma.$
Let $u\in U$ and suppose $\alpha(u)=\beta(x)$ for some $x\in \Gamma$. Claim $x\in U.$
If, instead, $x\notin U$ ($x\in\overline{\Gamma_0}$) then $x=\lim x_n$ for some sequence $\{x_n\}_n\subset \Gamma_0$ 
and each $x_n=g_{n1}\cdots g_{nk_{n}}$ for some $g_{nj}\in \fG_{m_j}$.
It may be assumed that $m_j\ne0$, since $g_{nj}$ would then be $1_\Gamma$ and hence, redundant. If $m_j>0$, then $g_{nj}\in\beta^{-1}(\alpha(\fG_{m_j-1}))$ and there is a $y_{nj}\in \fG_{m_j -1}$ such that $\beta(g_{nj})=\alpha(y_{nj})$ and and $\lim_n y_{n1}\cdots y_{nk_n}$ converges. If $m_j<0$, then there exists a $y_{nj}\in \fG_{m_j-1}=\alpha^{-1}(\beta(\fG_{m_j}))$ such that $\alpha(y_{nj})=\beta(g_{nj})$ and $\lim_n y_{n1}\cdots y_{nk_n}$ converges. Note if no such set of $y_{nj}$'s existed, then
no such set of $g_{nj}$'s could exist, but the $g_{nj}$'s do exist.
Therefore,
\begin{align*}
\beta(x)&=\beta(\lim x_n)=\lim \beta(x_n)\\
&=\lim \beta(g_{n1}\cdots g_{nk_n})\\
&=\lim\alpha(y_{n1}\cdots y_{nk_n})\\
&=\alpha(\lim y_{n1}\cdots y_{nk_n})\in\alpha(\overline{\Gamma_0}).
\end{align*}
That is, $\alpha(u)=\alpha(g)$ for some $g\notin U$; Hence, $ug^{-1}\in\ker\alpha\subseteq\Gamma_0$. This cannot be since then
$u=ug^{-1}g\in \Gamma_0\subseteq\overline{\Gamma_0}$. Hence, $x\in U$; that is, $\beta^{-1}(\alpha(u))\subseteq U$ whenever $u\in U$.

Now assume $\beta^{-1}(\alpha(x))\subseteq U$ for some $x\in \Gamma$. Then $\alpha(x)=\beta(u)$ for some $u\in U.$ A similar argument (with $\alpha$ and $\beta$ switched) shows that $x\in U.$ Hence, it has been shown that $U$ is a saturated hereditary open subset of $Q_{\alpha,\beta}(\Gamma).$\\\sq
\end{theorem}

\begin{corollary}\rm If $\alpha$ and $\beta$ are both automorphisms and $\Gamma$ is not the trivial group, 
then $\cO_{\alpha,\beta}(\Gamma)$ is not simple.\\
\pf Proceed by calculating $\Gamma_0$. Since $\ker\beta=\ker\alpha=\{1_\Gamma\}$, $\fG_n=\{1_\Gamma\}$ for all $n\in\zZ$. 
Thus, $\Gamma_0=\{1_\Gamma\},$ and is clearly not dense. Hence,
$Q_{\alpha,\beta}(\Gamma)$ is not minimal and by Theorem \ref{MuhlyThm}, $\cO_{\alpha,\beta}(\Gamma)$ is not simple.\\\sq
\end{corollary}

\begin{remark}\rm If $\Gamma$ is the trivial group, then the only endomorphism on $\Gamma$ is the identity automorphism, 
denoted $1$. Thus, $\Omega_{1,1}(\Gamma)=\{(1_\Gamma,1_\Gamma)\}$ and the topological group relation 
$Q_{1,1}(\Gamma)$ is the directed graph
$$\includegraphics[height=1cm]{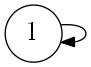}$$
Hence, $\cO_{1,1}(\Gamma)\cong\cO_1=C(\zT)$
which is not simple even though $\Gamma_0=\Gamma$  (it does not satisfy condition ($L$).)
\end{remark} 

Combining Theorem \ref{MinThm}, Theorem \ref{L} and Theorem \ref{MuhlyThm}, we obtain:

\begin{theorem}\label{SimpCrit}\rm \index{$\cO_{\alpha,\beta}(\Gamma)$! Simple}
If $\Gamma$ is a compact group with surjective endomorphisms $\alpha,\beta$ satisfying the conditions
\begin{enumerate}
\item $\abs{\ker\alpha},\abs{\ker\beta}<\infty$ and
\item either $\alpha$ or $\beta$ is not injective,
\end{enumerate}
then $\cO_{\alpha,\beta}(\Gamma)$ is simple if and only if $\Gamma_0$ is dense in $\Gamma.$
\end{theorem}

\begin{corollary}\rm If $\Gamma$ is a compact group with surjective endomorphisms $\alpha,\beta$ satisfying the conditions
\begin{enumerate}
\item $\abs{\ker\alpha},\abs{\ker\beta}<\infty$ and
\item either $\alpha$ or $\beta$ is not injective,
\end{enumerate}
then $\cO_{\alpha,\beta}(\Gamma)$ is simple if and only if $\cO_{\beta,\alpha}(\Gamma)$ is.
\end{corollary}

\begin{definition}\cite{RS} \rm A \emph{Kirchberg algebra}\index{Kirchberg Algebra} 
is a purely infinite, simple, nuclear, separable $C^*$-algebra.
\end{definition}

\begin{corollary}\label{OFG} \rm Let $F,G\in M_d(\zZ)$ have non-zero determinants and assume either $\abs{\det F}\ne 1$ 
or $\abs{\det G}\ne1.$  
Then if $\Gamma_0\subseteq\zT^d$ is dense, then $\cO_{F,G}(\zT^d)$ is a unital classifiable Kirchberg algebra.\\
\pf By Theorem \ref{SimpCrit}, $\cO_{F,G}(\zT^d)$ is simple. Furthermore, $\cO_{F,G}(\zT^d)$ is unital and separable.
By \cite[Proposition 5.2, Theorem 5.6, and Theorem 5.9]{Mc1}, we have
$\cO_{F,G}(\zT^d)$ is nuclear, (and along with \cite[Page 76]{RS} and \cite[Proposition 3.3]{S}) purely infinite, 
and in the UCT class. Hence, $\cO_{F,G}(\zT^d)$ is a classifiable Kirchberg algebra.\\\sq 
\index{$\cO_{F,G}(\zT^d)$! Purely Infinite}\index{$\cO_{F,G}(\zT^d)$! Simple}
\end{corollary}

\begin{proposition}\cite[Proposition 4.2]{EaHR}\label{Raeburn} \rm Let $F$ be an integer dilation matrix; that is, an integral matrix
whose eigenvalues have modulus greater than 1. Then $\cO_{F,1}(\zT^d)$ is simple and purely infinite.\index{Integer Dilation Matrix}
\end{proposition}

\begin{corollary}\rm For integer dilation matrix $F$, we have $\Gamma_0$ is dense in $\zT^d$ and so
$\cO_{F,1}(\zT^d)$ and $\cO_{1,F}(\zT^d)$ are unital classifiable Kirchberg algebras.
\end{corollary}

It becomes difficult to calculate the density of $\Gamma_0$ in practice, but the next few theorems will aid this calculation immensely.

\begin{theorem}\rm Let $\Gamma$ be a compact group and let $\alpha,\beta,\gamma\in\mbox{End}(\Gamma)$ be surjective.
Then
\begin{enumerate}[(a)]
\item $\cO_{\alpha,\beta}(\Gamma)$ is simple if and only if $\cO_{\alpha\circ\gamma,\beta\circ\gamma}(\Gamma)$ is simple.
\item $\cO_{\alpha,\beta}(\Gamma)$ is simple guarantees $\cO_{\gamma\circ\alpha,\gamma\circ\beta}(\Gamma)$ is simple.
\end{enumerate}
\pf (a) Let 
$$\fG_n=\begin{cases}\alpha^{-1}(\beta(\fG_{n-1}))&\mbox{if }n>0\\\beta^{-1}(\alpha(\fG_{n+1}))&\mbox{if }n<0\\
\{1_\Gamma\}&\mbox{if } n=0\end{cases}$$
and let 
$$\fG_n'=\begin{cases}(\alpha\circ\gamma)^{-1}(\beta\circ\gamma)(\fG_{n-1}')&\mbox{if }n>0\\
(\beta\circ\gamma)^{-1}(\alpha\circ\gamma)(\fG_{n+1}')&\mbox{if }n<0\\
\{1_\Gamma\}&\mbox{if } n=0.\end{cases}$$
Then
$$\fG_1'=(\alpha\circ\gamma)^{-1}(\beta\circ\gamma)(1_\Gamma)=(\alpha\circ\gamma)^{-1}(1_\Gamma)
=\gamma^{-1}(\ker\alpha)=\gamma^{-1}(\fG_1).$$
Next, assume $\fG_n'=\gamma^{-1}(\fG_n)$ for all $n\le k-1.$ Then
\begin{align*}
\fG_k'&=(\alpha\circ\gamma)^{-1}(\beta\circ\gamma)(\fG_{k-1})\\
&=(\alpha\circ\gamma)^{-1}(\beta\circ\gamma)(\gamma^{-1}(\fG_{k-1})\\
&=(\alpha\circ\gamma)^{-1}(\beta(\fG_{k-1}))\\
&=\gamma^{-1}(\alpha^{-1}(\beta(\fG_{k-1})))\\
&=\gamma^{-1}(\fG_k).
\end{align*}
Hence, by induction, $\fG_n'=\gamma^{-1}(\fG_n)$. Note one can use the same argument with $\alpha$ and $\beta$ flipped 
for $n<0.$
Next, $\subgp{\gamma^{-1}(\cup_n\fG_n)}=\gamma^{-1}(\subgp{\cup_n\fG_n}).$ To see this claim, 
if $x\in\subgp{\gamma^{-1}(\cup_n\fG_n)}$ then $x=\prod_{m=1}^k x_m$ where $x_m\in\gamma^{-1}(\fG_{n(m)}).$
That is, $\gamma(x_m)\in\fG_{n(m)}$ so
$$\gamma(x)=\prod_{m=1}^k \gamma(x_m)\in\subgp{\cup_n\fG_n}.$$
Conversely, if $\gamma(x)\in\subgp{\cup_n\fG_n}$ then $\gamma(x)=\prod_{m=1}^k y_m$ such that $y_m\in\fG_{n(m)}$ for each $m=1,...,k.$ Note that $\prod_{m=1}^k\gamma^{-1}(y_m)=(\prod_{m=1}^ky_m')\ker\gamma$ for some $y_m'\in\gamma^{-1}(y_m)\subset\fG_{n(m)}'.$ If $y=\prod_{m=1}^k y_m',$ then $y^{-1}x\in\ker\gamma$ and hence,
$$x=y\cdot y^{-1}x\in y\cdot\ker\gamma=\prod_{m=1}^k\gamma^{-1}(y_m)\subset\subgp{\gamma^{-1}(\cup_n\fG_n)}.$$
Thus, $$\subgp{\gamma^{-1}(\cup_n\fG_n)}=\gamma^{-1}(\subgp{\cup_n\fG_n}).$$

Next, $\overline{\gamma^{-1}(\subgp{\cup_n\fG_n})}=\gamma^{-1}(\overline{\subgp{\cup_n\fG_n}}).$ To see this claim,
if $x\in \overline{\gamma^{-1}(\subgp{\cup_n\fG_n})}$ then $x=\lim_m x_m$ where $x_m\in\gamma^{-1}(\subgp{\cup_n\fG_n}).$ Then
$$\gamma(x)=\lim_m\gamma(x_m)\in\overline{\subgp{\cup_n\fG_n}}.$$
Conversely, assume $\gamma(x)\in\overline{\subgp{\cup_n\fG_n}}$ but $x\notin\overline{\gamma^{-1}(\subgp{\cup_n\fG_n})}.$ Then $\gamma(x)=\lim_m x_m$ for some $x_m\in \subgp{\cup_n\fG_n}$ and there exists an
open set $U\subset \Gamma$ such that $x\in U$ and
$$U\cap \gamma^{-1}(\subgp{\cup_n\fG_n})=\emptyset.$$
If $z\in \gamma(U)\cap\subgp{\cup_n\fG_n}$ then $z=\gamma(u)$ for some $u\in U$ and so $u\in\gamma^{-1}(\subgp{\cup_n\fG_n}).$ This is impossible, so
$$\gamma(U)\cap\subgp{\cup_n\fG_n}=\emptyset.$$
But $\gamma(x)\in\overline{\subgp{\cup_n\fG_n}}$ and $\gamma(U)$ is open imply $\gamma(x)\notin\gamma(U).$ This
is a contradiction since $x\in U.$ Thus, it must be that
$$x\in\overline{\gamma^{-1}(\subgp{\cup_n\fG_n})}$$
and so
$$\overline{\gamma^{-1}(\subgp{\cup_n\fG_n})}=\gamma^{-1}(\overline{\subgp{\cup_n\fG_n}}).$$

Finally, 
\begin{align*}
\cO_{\alpha,\beta}(\Gamma)\mbox{ is simple }&\mbox{if and only if }\overline{\subgp{\cup_n\fG_n}}=\Gamma\\
&\mbox{if and only if }\Gamma=\gamma^{-1}(\Gamma)=\gamma^{-1}(\overline{\subgp{\cup_n\fG_n}})
=\overline{\gamma^{-1}(\subgp{\cup_n\fG_n})}\\
&\hspace{1.07in}=\overline{\subgp{\gamma^{-1}(\cup_n\fG_n)}}=\overline{\subgp{\cup_n\gamma^{-1}(\fG_n)}}=\overline{\subgp{\cup_n\fG_n'}}\\
&\mbox{if and only if }\cO_{\alpha\circ\gamma,\beta\circ\gamma}(\Gamma)\mbox{ is simple.}
\end{align*}

(b) Let $$\fG_n'=\begin{cases}(\gamma\circ\alpha)^{-1}(\gamma\circ\beta)(\fG_{n-1}')&\mbox{if }n>0\\
(\gamma\circ\beta)^{-1}(\gamma\circ\alpha)(\fG_{n+1}')&\mbox{if }n<0\\
\{1_\Gamma\}&\mbox{if } n=0.\end{cases}$$
Note $\fG_1=\ker\alpha\subset \ker(\gamma\circ\alpha)=\fG_1'.$ Now assume $\fG_n\subset\fG_n'$ for all $0<n\le k-1.$ Then
if $x\in\fG_k$ then there exists $y\in\fG_{k-1}\subset\fG_{k-1}'$ such that $\alpha(x)=\beta(y).$ Thus,
$(\gamma\circ\alpha)(x)=(\gamma\circ\beta)(y)$ and so
$$x\in (\gamma\circ\alpha)^{-1}(\gamma\circ\beta)(y)\subset (\gamma\circ\alpha)^{-1}(\gamma\circ\beta)(\fG_{k-1}')=\fG_k'.$$ 
Hence, by induction, $\fG_n\subset\fG_n'$. The argument for $n<0$ is identical. Thus, since $\subgp{\cup_n\fG_n}$ is dense
in $\Gamma$ and $\subgp{\cup_n\fG_n}\subset\subgp{\cup_n\fG_n'},$ it must be that
$\subgp{\cup_n\fG_n'}$ is dense in $\Gamma$ and hence,
$$\cO_{\gamma\circ\alpha,\gamma\circ\beta}(\Gamma)\mbox{ is simple}.$$
\sq
\end{theorem} 

\begin{corollary}\label{SimpleFG}\rm Given $d\times d$-matrices $F$ and $G$ with non-zero determinant,
$$\cO_{F,G}(\zT^d)\mbox{ is simple if and only if }\cO_{\det F,G(\mbox{\tiny adj }F)}(\zT^d)\mbox{ is simple}$$
and likewise,
$$\cO_{F,G}(\zT^d)\mbox{ is simple if and only if }\cO_{\det F,(\mbox{\tiny adj }F)G}(\zT^d)\mbox{ is simple}$$
where $\mbox{adj }F$ is the adjugate of $F.$\\
\pf Let $H=\mbox{adj }F$ and $H'=\mbox{adj G}$ then
\begin{align*}
\cO_{F,G}(\zT^d)\mbox{ is simple }&\implies \cO_{HF,HG}(\zT^d)\mbox{ is simple }&(\mbox{left multiply by }H)\\
&\iff\cO_{\det F, HG}(\zT^d)\mbox{ is simple }&\\
&\iff\cO_{(\det F)H', (\det G)H}(\zT^d)\mbox{ is simple }&(\mbox{right multiply by }H')\\
&\iff\cO_{(\det F)H'F,(\det G)(\det F)}(\zT^d)\mbox{ is simple }&(\mbox{right multiply by }F)\\
&\iff\cO_{H'F,\det G}(\zT^d)\mbox{ is simple }&(\mbox{right factor by }\det F)\\
&\implies\cO_{(\det G)F,(\det G )G}(\zT^d)\mbox{ is simple }&(\mbox{left multiply by }G)\\
&\iff\cO_{F,G}(\zT^d)\mbox{ is simple }&(\mbox{right factor by }\det G)\\
\end{align*}
Hence, 
$$\cO_{F,G}(\zT^d)\mbox{ is simple if and only if }\cO_{\det F, HG}(\zT^d)\mbox{ is simple }$$
and likewise,
$$\cO_{F,G}(\zT^d)\mbox{ is simple if and only if }\cO_{\det F, GH}(\zT^d)\mbox{ is simple }$$\sq
\end{corollary}

\begin{corollary}\rm Given $d\times d$-matrices $F,$ $G,$ and $H$ with non-zero determinant,
$$\cO_{F,G}(\zT^d)\mbox{ is simple if and only if }\cO_{HF,HG}(\zT^d)\mbox{ is simple.}$$
\pf We proceed similar to the last corollary.
\begin{align*}
\cO_{F,G}(\zT^d)\mbox{ is simple }&\implies \cO_{HF,HG}(\zT^d)\mbox{ is simple }&(\mbox{left multiply by }H)\\
&\iff \cO_{\det H\det F,\det H(\mbox{\tiny adj }F)G}(\zT^d)\mbox{ is simple }& (\mbox{by Corollary \ref{SimpleFG}})\\
&\iff \cO_{\det F,(\mbox{\tiny adj }F)G}(\zT^d)\mbox{ is simple }&(\mbox{right factor by }\det H)\\
&\iff \cO_{F,G}(\zT^d)\mbox{ is simple }&(\mbox{by Corollary \ref{SimpleFG}})\\
\end{align*}\sq
\end{corollary}

\begin{corollary}\rm Given $d\times d$-matrices $F$ and $G$ with non-zero determinants, there exists unimodular matrices $U$ and
$V$ and a positive diagonal matrix $D$ such that
$$(\mbox{adj }F)G=UDV$$
Therefore,
$$\cO_{F,G}(\zT^d)\mbox{ is simple if and only if }\cO_{\det F, DVU}(\zT^d)\mbox{ is simple.}$$
\pf By Corollary \ref{SimpleFG},
$$\cO_{F,G}(\zT^d)\mbox{ is simple if and only if }\cO_{\det F,(\mbox{\tiny adj }F)G}(\zT^d)\mbox{ is simple}$$
and by \cite[Proposition 3.19]{Mc1},
$$\cO_{\det F,UDV}(\zT^d)\cong\cO_{U^{-1}\det F,DV}(\zT^d)\cong\cO_{\det F,DVU}(\zT^d)$$\sq
\end{corollary}

\begin{remark}\rm The previous corollary states that one need only characterize the case where $F$ is a multiple of the identity
and $G$ is a positive diagonal matrix multiplied by a unimodular matrix.
\end{remark} 

\begin{lemma}\rm If $x\in\overline{\Gamma_0}$ and $y\in\fG_n,$ then $xy\in\overline{\Gamma_0}.$\\
\pf Let $x=\lim_k x_k\in\overline{\Gamma_0}$ where the $x_k$'s are products of elements in $\{\fG_m\}_{m\in\zZ}.$
Then for any $y\in\fG_n,$ $x_ky$ is a product of elements in $\{\fG_m\}_{m\in\zZ}$ that converge to $xy.$
That is,
$$xy=\lim_k x_ky\in\overline{\Gamma_0}.$$
\sq
\end{lemma}

\begin{lemma}\rm\label{G0} Assume $\Gamma=\Gamma_1\times\Gamma_2$ is a compact group. Furthermore,
assume $\Delta_i\subset\Gamma_i$  for $i=1,2$ have the following properties:
\begin{enumerate}
\item $\Delta_1\times \{1_{\Gamma_2}\}\subset \overline{\Gamma_0},$
\item for every $z\in\Delta_2,$ there is a $w\in\Delta_1$ such that $(w,z)\in\fG_n$ for some $n\in\zZ.$
\end{enumerate}
Then
$$\{1_{\Gamma_1}\}\times \Delta_2\subset\overline{\Gamma_0}$$
and hence,
$$\Delta_1\times\Delta_2\subset\overline{\Gamma_0}.$$
\pf Let $z\in\Delta_2.$ Then there exists $w\in\Delta_1$ such that $(w,1_{\Gamma_2})\in\overline{\Gamma_0}$ and $(w,z)\in\fG_n$
for some $n\in\zZ.$ Thus, by the previous lemma,
$$(1_{\Gamma_1},z)=(ww^{-1},z)=(w,z)(w,1_{\Gamma_2})^{-1}\in\overline{\Gamma_0}$$ 
since $(w,1_{\Gamma_2})^{-1}\in\overline{\Gamma_0}$ and $\Gamma_0$ is a group. As $z\in\Delta_2$ was arbitrary,
$$\{1_{\Gamma_1}\}\times\Delta_2\subset \overline{\Gamma_0}$$
and hence,
$$\Delta_1\times \Delta_2\subset\overline{\Gamma_0}.$$\sq
\end{lemma}

By Corollary \ref{SimpleFG}, one can reduce the study of simplicity of $\cO_{F,G}(\zT^d)$ to that of the case 
$\cO_{n,G}(\zT^d)$ where $n$ is a positive integer. In turn, we now consider the case where $G$ is upper (or lower) triangular, 
but we shall need the following lemma.

\begin{lemma}\rm If $a,b\in\zN$ such that $a>b$ and $k=\gcd(a,b),$ then there exists $c\in\zZ$ such that
$$(\frac{b}{k})c=1\mod(\frac{a}{k})$$
\pf Since $\gcd(\frac{a}{k},\frac{b}{k})=1,$ there exist $m,c\in\zZ$ such that
$$(\frac{b}{k})c+(\frac{a}{k})m=1$$
Hence,
$$(\frac{b}{k})c=1\mod(\frac{a}{k})$$\sq
\end{lemma}

\begin{theorem}\rm Suppose $G=(b_{ij})_{i,j=1}^d$ is an upper (or lower) triangular $d\times d$-matrix with non-zero determinant.  If $\abs{b_{jj}}\ne n$ for each $j=1,...,d,$ then $\cO_{n,G}(\zT^d)$ is simple.\\
\pf We shall only prove the upper triangular case (the lower triangular case is similar) and we shall view the torus $\zT$ as $\zR/\zZ$
for the convenience of this calculation. To this end, let $k_j=\gcd(n,\abs{b_{jj}})$
and let $e_j$ denote the standard vector basis with a 1 in the $j$-th coordinate and 0 elsewhere.
First suppose $n>\abs{b_{11}}.$ Then $\frac{1}{n}e_1\in\fG_1$ and $\frac{b_{11}}{n}e_1\in G(\fG_1).$ By the previous lemma,
there exists $c_{11}\in\zZ$ such that 
$$\frac{b_{11}}{k_1}c_{11}=1\mod(\frac{n}{k_1})$$
and hence,
$$\frac{1}{n/k_1}e_1\in G(\fG_1), \frac{1}{n(n/k_1)}e_1\in\fG_2,\mbox{ and } \frac{b_{11}}{n(n/k_1)}e_1\in G(\fG_2).$$
Then by the previous lemma again, there exists $c_{12}\in\zZ$ such that 
$$\frac{b_{11}}{n(n/k_1)}c_{12}=1\mod(\frac{n^2}{k_1^2})$$
and hence,
$$\frac{1}{n^2/k_1^2}e_1\in G(\fG_2)\mbox{ and }\frac{1}{n(n^2/k_1^2)}e_1\in\fG_3.$$
Continue by induction to obtain
$$\frac{1}{n^m/k_1^m}e_1\in G(\fG_m)\mbox{ and }\frac{1}{n(n^m/k_1^m)}e_1\in\fG_{m+1}$$
implying 
$$\zT\times 0^{d-1}\subseteq\overline{\Gamma_0}.$$
If $n<\abs{b_{11}}$ then $\frac{1}{b_{11}}e_1\in\fG_{-1}$ and $\frac{n}{b_{11}}\in F(\fG_{-1}).$ Then by the previous lemma,
there exists $c_{11}\in\zZ$ such that $$(\frac{n}{k_1})c_{11}=1\mod(\frac{b_{11}}{k_1})$$
and hence,
$$\frac{1}{b_{11}/k_1}e_1\in F(\fG_{-1})$$
and proceed as before to obtain
$$\frac{1}{b_{11}^m/k_1^m}e_1\in F(\fG_{-m})\mbox{ and } \frac{1}{b_{11}(b_{11}^m/k_1^m)}e_1\in\fG_{-m-1}$$
and hence,
$$\zT\times 0^{d-1}\subseteq \overline{\Gamma_0}.$$
In either case,
$$\zT\times 0^{d-1}\subseteq \overline{\Gamma_0}.$$

Next, if $n>\abs{b_{22}}$ then $\frac{1}{n}e_2\in\fG_1$ and $\frac{b_{12}}{n}e_1+\frac{b_{22}}{n}e_2\in G(\fG_1).$
Then
$$\frac{b_{12}}{n^2}e_1+\frac{b_{22}}{n^2}e_2\in \fG_2$$
and since $\zT\times 0^{d-1}\subseteq\overline{\Gamma_0},$ Lemma \ref{G0} guarantees 
$$\frac{b_{22}}{n^2}e_2\in\overline{\Gamma_0}$$
and by the previous lemma, there exists $c_{21}\in\zZ$ such that
$$(\frac{b_{22}}{k_2})c_{21}=1\mod(\frac{n^2}{k_2})$$
hence,
$$\frac{1}{n^2/k_2}e_2\in\overline{\Gamma_0}.$$
Thus, 
$$\frac{b_{12}}{n^2/k_2}e_1+\frac{b_{22}}{n^2/k_2}e_2\in G(\overline{\Gamma_0})$$
and
$$\frac{b_{12}}{n^3/k_2}e_1+\frac{b_{22}}{n^3/k_2}e_2\in \overline{\Gamma_0}$$
Again Lemma \ref{G0} guarantees that
$$\frac{b_{22}}{n^3/k_2}e_2\in\overline{\Gamma_0}$$
and the previous lemma may be used again to obtain
$$\frac{1}{n^3/k_2^2}e_2\in\overline{\Gamma_0}.$$
One may proceed by induction to obtain
$$\frac{1}{n^{m+1}/k_2^m}e_2\in\overline{\Gamma_0}$$
and hence,
$$0\times\zT\times 0^{d-2}\subseteq\overline{\Gamma_0}.$$
If $n<\abs{b_{22}}$ then a similar argument shows $$0\times\zT\times 0^{d-2}\subseteq\overline{\Gamma_0}$$
and in either case,
$$0\times\zT\times 0^{d-2}\subseteq\overline{\Gamma_0}.$$
Proceed with $b_{33},...,b_{dd}$ to obtain
$$\zT\times 0^{d-1},\, 0\times\zT\times 0^{d-2},\,...,\,0^{d-1}\times\zT\subseteq\overline{\Gamma_0}$$
and hence,
$$\overline{\Gamma_0}=\zT^d$$\sq
\end{theorem}

\begin{theorem}\rm Let $\Gamma=\zT^d\cong\zR^d/\zZ^d.$ If there is an open set $U\subset\overline{\Gamma_0},$
then $$\overline{\Gamma_0}=\zT^d.$$
\pf Suppose $U\subset\overline{\Gamma_0}$ is open. Then there exist $\bar{N}=(N_j)_{j=1}^d\in\zZ^d$ and
$$x_{\bar{n}}=\sum_{j=1}^d (\frac{1}{2^{n_j}}+\frac{k_{n_j,2}}{2^{n_j-1}}+...+\frac{k_{n_j,n_j-1}}{2})e_j\in U$$
for all $\bar{n}=(n_j)_{j=1}^d\ge\bar{N}$
where $k_{n_j,i}\in\{0,1\}$ for each $i,j$ and $e_j=(0,...,0,1,0,...,0)$ with a $1$ in the $j$-th coordinate. That is,
there exists $x\in U$ such that, for all $j=1,...,d$
$$\{x+(\frac{k_{N_j}}{2^{N_j+1}}+...+\frac{k_{N_j+m}}{2^{N_j+m}})e_j\st m\in\zN, k_{N_j+p}\in\{0,1\}\mbox{ for all }p=0,...,m\}\subset U.$$
Hence,
$$\frac{1}{2^{N_j+m}}e_j=(x+\frac{1}{2^{N_j+m}}e_j)-x\in\subgp{U}\subset\overline{\Gamma_0}$$
for each $m\in\zN.$ 
Since $\{\frac{1}{2^{N_j+m}}\}_{m\in\zN}$ is dense in $\zT$, we obtain
$$0^{j-1}\times \zT\times 0^{d-j}\subset\overline{\Gamma_0}$$
for each $j=1,...,d$ and hence
$$\overline{\Gamma_0}=\zT^d.$$\sq
\end{theorem}

\section{Acknowledgements}This paper was based upon my thesis work at the University of Calgary and research done during a post-doctoral fellowship at the University of Regina. I am indebted to NSERC, the Department of Mathematics and Statistics at the
University of Calgary and the Department of Mathematics and Statistics at the University of Regina 
in providing funding to finance my mathematical studies.

%%%%%%%%%%%%%%%%%%%%%%%%%%%%%%%%%

\end{document}